\documentclass[12pt]{article}
\usepackage{fancyhdr,amsmath, graphicx, psfrag, color ,amsfonts, layout}

\font\fontauthors=cmcsc10 scaled \magstep1

 at 8truept
 at 8truept

\newtheorem{Th}{Theorem}

\newtheorem{Lem}[Th]{Lemma}

\newtheorem{Rem}[Th]{Remark}

\newcommand{\pff}{\noindent {\sc Proof.}\ }

\def\QED{\hfill$\square$ \vskip 10pt}

\newtheorem{prerem}[Th]{Remark}

\def\be{\begin{equation}}

\def\ee{\end{equation}}
\def\egalLoi{{~\mathop{= }\limits^{\rond L}}~}

\def\Supp{\mathop{\rm Supp}\nolimits}
\def\stirling2 #1#2{\left\{\begin{matrix} #1\\#2\end{matrix}\right\}}

\def\Var{\mathop{\rm Var}\nolimits}

\newcommand{\C}{\;{}^{{}_\vert}\!\!\!{\rm C}}

\newcommand{\Q}{{\rm Q}\kern-.65em{}^{{}_/}}

\def\1{{\bf 1}}

\def\g#1{\mathbb #1}
\def\rond#1{\mathcal #1}

\def\C{\g C}

\begin{document}

\begin{center}
\LARGE{\bf {Support and density of the limit  $m$-ary search trees distribution

}}
\end{center}

\begin{center}
{\fontauthors
Brigitte Chauvin}\footnote{
Universit\'e de Versailles-St-Quentin,
Laboratoire de Math\'ematiques de Versailles,
CNRS, UMR 8100,
45, avenue des Etats-Unis, 78035 Versailles CEDEX, France.
}

\medskip
{\fontauthors
Quansheng Liu}\footnote{
LMAM,
Universit\'e  de Bretagne Sud,
Campus de Tohannic, BP 573, 56017 Vannes, France.
}

\medskip
{\fontauthors
Nicolas Pouyanne}\footnote{
Universit\'e de Versailles-St-Quentin,
Laboratoire de Math\'ematiques de Versailles,
CNRS, UMR 8100,
45, avenue des Etats-Unis, 78035 Versailles CEDEX, France.
}

\end{center}

\begin{center}
{\it 17 janvier 2012}\\
\end{center}
{\small
\noindent{\bf Abstract.}
The space requirements of an $m$-ary search tree satisfies a well-known phase transition: when $m\leq 26$, the second order asymptotics is Gaussian. When $m\geq 27$, it is not Gaussian any longer and a limit $W$ of a complex-valued martingale arises. We show that the distribution of $W$ has a square integrable density on the complex plane, that its support is the whole complex plane, and that it has finite exponential moments. The proofs are based on the study of the distributional equation $ W\egalLoi\sum _{k=1}^mV_k^{\lambda}W_k$, where $V_1, ..., V_m$ are the spacings of $(m-1)$ independent random variables uniformly distributed on $[0,1]$, $W_1, ..., W_m$ are independent copies of W which are also independent of $(V_1, ..., V_m)$ and $\lambda$ is a complex number.

\bigskip

Smoothing transformation. Absolute continuity. Support. Exponential moments. Mandelbrot cascades.

}


\section{Introduction}\label{intro}

Search trees are fundamental data
structures in computer science used in searching and sorting. For integers $m\geq 2$, $m$-ary search trees generalize the binary search tree. The
quantity $m$ is called the branching factor.

A random $m$-ary search tree is an $m$-ary tree in which each node has the capacity to
contain $(m-1)$ elements called the data or the keys. The keys can be considered as i.i.d. random variables $x_i,~i\geq 1$, with any diffusive distribution on the interval $[0, 1]$.

The tree $T_n, n\geq 0$, is recursively defined as follows: $T_0$ is reduced to an empty node-root; $T_1$ is reduced to a node-root which contains $x_1$, $T_2$ is reduced to a node-root which contains $x_1$ and $x_2$, \dots , $T_{m-1}$ still has one node-root, containing $x_1, \dots x_{m-1}$.   As soon as the ($m-1$)-th key is inserted in the root, $m$ empty subtrees of the root are created, corresponding from left to right to the $m$ ordered intervals $I_1 = ]0,x_{(1)}[, \dots, I_m=]x_{(m-1)}, 1[$, where $0< x_{(1)} < \dots < x_{(m-1)} <1$ are the ordered first $(m-1)$ keys. Each following key $x_m, \dots $ is recursively inserted in the subtree corresponding to the unique interval $I_j$ to which it belongs. As soon as a node is saturated, $m$ empty subtrees of this node are created.

For each $i=\{1,\dots ,m-1\}$ and $n\geq 1$, $X_n^{(i)}$ is the number
of nodes in $T_n$ which contain $(i-1)$ keys (and $i$ gaps or free places) after insertion of the $n$-th key; such nodes are named nodes of type $i$. We only take into consideration the external nodes and not the internal nodes which are the saturated nodes. The vector $X_n$ is called  the composition vector of the $m$-ary search tree. It provides a model for the space requirement of the algorithm. By spreading the input
data in $m$ directions instead of only $2$, as is the case for a binary search tree, one
seeks to have shorter path lengths and thus quicker searches. One can refer to Mahmoud's book \cite{Mahmoud} for further details on search trees.


\medskip
The following figure is an example of $4$-ary search tree obtained by insertion of the successive
numbers $0.3$ ,$0.1$, $0.4$, $0.15$, $0.9$, $0.2$, $0.6$, $0.5$, $0.35$, $0.8$, $0.97$, $0.93$, $0.23$, $0.84$, $0.62$, $0.64$, $0.33$, $0.83$.
The corresponding composition vector is $X_{18}=\!~^t(9,2,2)$.

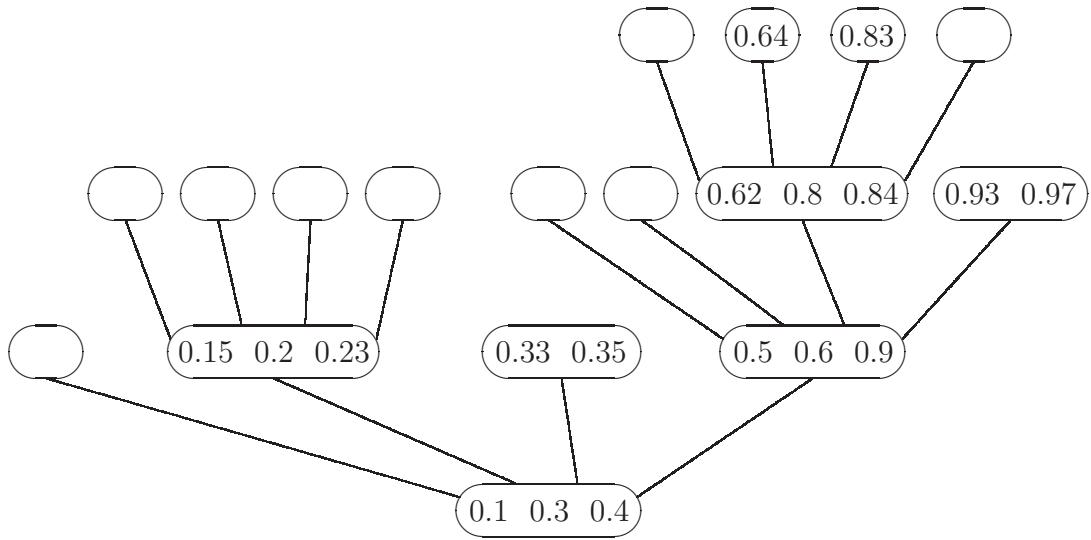
\begin{picture}(400,210)
\put(180,0){$0.1~~0.3~~0.4$}
\put(210,4){\oval(70,20)}
\put(20,64){\oval(28,20)}
\put(70,60){$0.15~~0.2~~0.23$}
\put(106,64){\oval(80,20)}
\put(190,60){$0.33~~0.35$}
\put(215,64){\oval(60,20)}
\put(280,60){$0.5~~0.6~~0.9$}
\put(310,64){\oval(70,20)}
\put(50,124){\oval(28,20)}
\put(85,124){\oval(28,20)}
\put(120,124){\oval(28,20)}
\put(155,124){\oval(28,20)}
\put(210,124){\oval(28,20)}
\put(245,124){\oval(28,20)}
\put(270,120){$0.62~~0.8~~0.84$}
\put(306,124){\oval(80,20)}
\put(360,120){$0.93~~0.97$}
\put(385,124){\oval(58,20)}
\put(251,184){\oval(28,20)}
\put(280,180){$0.64$}
\put(291,184){\oval(28,20)}
\put(320,180){$0.83$}
\put(331,184){\oval(28,20)}
\put(371,184){\oval(28,20)}
\qbezier(176.5,9)(98.75,31.5)(20,54)
\qbezier(198,14)(152,34)(106,54)
\qbezier(221,14)(218,34)(215,54)
\qbezier(243.5,9)(276.75,31.5)(310,54)
\qbezier(67,69)(58.5,91.5)(50,114)
\qbezier(94,74)(89.5,94)(85,114)
\qbezier(118,74)(119,94)(120,114)
\qbezier(145,69)(150,91.5)(155,114)
\qbezier(276,69)(243,91.5)(210,114)
\qbezier(299,74)(273,94)(245,114)
\qbezier(322,74)(314,94)(306,114)
\qbezier(344,69)(364.5,91.5)(385,114)
\qbezier(267,129)(259,151.5)(251,174)
\qbezier(295,134)(293,154)(291,174)
\qbezier(317,134)(324,154)(331,174)
\qbezier(345,129)(358,151.5)(371,174)
\end{picture}

\vskip 20pt
A numerous literature is devoted to the asymptotic behavior of this composition vector.
A famous phase transition appears. When $m\leq 26$, the random vector admits a central limit theorem
with convergence in distribution to a Gaussian vector: see also Mahmoud's book \cite{Mahmoud} or Janson \cite{Jan} for a vectorial treatment.
\medskip

When  $m\geq 27$, an almost sure asymptotics for the composition vector has been obtained in \cite{ChaPou04}:
\begin{equation}
\label{discrete}
X_n= nv_1 +  \Re(n^{\lambda_2} W v_2 ) + o(n^{\sigma_2}) \quad a.s.,
\end{equation}
where $\lambda_2= \sigma_2 + i \tau_2 $ is the root of the polynomial
\be
\label{car}
\prod _{k=1}^{m-1}(z +k)-m!
\ee
having the second largest real part
$\sigma_2$
and a positive imaginary part $\tau_2$,
$v_1$ and $v_2$ are two deterministic vectors,  and $W$ is the limit of a complex-valued martingale that admits moments of all positive orders.

Heated conjectures  about the  second order complex-valued limit distribution $W$ remain open
(see \cite{ChaPou04}, \cite{Pou05}, Chern and Hwang \cite{ChernHwang}, Mahmoud \cite{Mahmoud},  Janson \cite{Jan}).

\medskip
A significant step is achieved by Fill and Kapur in \cite{FillKapur}, who establish that $W$ satisfies the following distributional equation called the {\em smoothing equation}:
\begin{equation}
\label{distribW}
W\egalLoi\sum _{k=1}^mV_k^{\lambda_2}W_k,
\end{equation}
where $V_1, ..., V_m$ are the spacings of $(m-1)$ independent random variables uniformly distributed on $[0,1]$,
 $W_1, ..., W_m$ are independent copies of $W$ which are also independent of $(V_1, ..., V_m)$.
The precise definition  of $V_j$ will be given  hereunder.
By a contraction method, Fill and Kapur prove that $W$ is the unique solution of Eq.~(\ref{distribW}) in the space $\rond M_2(C)$ of square integrable probability measures having
$C=\g E(W)$ as expectation.
The present paper is based on this characterization of $W$.

\medskip

It has been recently proved \cite{ChaLiuPou} that the continuous-time embedding of the process $(X_n)_n$ has an analogous asymptotic behavior, with a second-order term which is a solution of some distributional equation (not the same one). Inspired by this study of the continuous-time case 
we prove the following theorem.

\begin{Th}
\label{W}
Let $W$ be the second order limit distribution of an $m$-ary search tree for $m\geq 27$, defined by (\ref{discrete}).

\vskip 5pt
\hskip 10pt (i) The support of $W$ is the whole complex plane.

\vskip 5pt
\hskip 10pt (ii) The law of $W$ admits a continuous square integrable density on $\g C$.

\vskip 5pt
\hskip 10pt (iii)
 $\g E e^{\delta  |W|} < \infty$ for some $\delta >0$.
The exponential moment generating series of~$W$ (thus) has a positive radius of convergence.

\end{Th}
Thanks to Fill and Kapur results \cite{FillKapur}, these results are immediate corollaries of Theorems \ref{support} and  \ref{density} proved in the next two sections.

\medskip
In the whole sequel, let $V_1,\dots ,V_m$ be the spacings of $(m-1)$ independent  random variables uniformly distributed on $[0,1]$.
In other words, let $U_1, \dots , U_{m-1}$ be independent  random variables uniformly distributed on  $[0,1]$ and let
$U_{(1)}\leq \dots \leq U_{(m-1)}$ be their order statistics.
Denote also $U_{(0)} := 0$, $U_{(m)} :=1$.
For any $k\in\{ 1,\dots ,m\}$, the random variable $V_k$ is defined by
$$
V_k:= U_{(k)} - U_{(k-1)}.
$$
The variables $V_k$ are $Beta(1,m-1)$-distributed and satisfy $\sum_{k=1}^{m} V_k=1$ almost surely.
\begin{Rem}
Details about roots of (\ref{car}) can be found in Hennequin \cite{HennequinThese} and Mahmoud \cite{Mahmoud}.
Note that for $m=2$, the polynomial   (\ref{car}) has the unique root $\lambda = 1$. For $m\geq 3$, it is known that if $\lambda_2$ is a root of the polynomial (\ref{car})
 having the second largest real part, 
then 
 $\lambda_2$ is non real, $\Re\lambda_2 <1$ for any $m\geq 3$, $\Re\lambda_2 > 0$ if and only if $m\geq 14$, and
$$
\Re\lambda_2 >\frac 12\Longleftrightarrow m\geq 27.
$$
\end{Rem}

\section{Support}
\label{sec:support}

The limit distribution $W$ satisfies (\ref{distribW}). From now on, we consider the solutions of the distributional equation
\begin{equation}
\label{smoothing}
Z\egalLoi\sum _{k=1}^mV_k^{\lambda}Z_k,
\end{equation}
where $V_1, ..., V_m$ are the spacings of $(m-1)$ independent random variables uniformly distributed on $[0,1]$, $Z_1, ..., Z_m$ are independent copies of $Z$ which are also independent of $(V_1, ..., V_m)$ and $\lambda$ is a non real complex number.

\medskip

We assume that
\be
\label{hyp}
\lambda \hbox{ \it is a non real root  of (\ref{car}) having a positive real part } \sigma.
\ee

Indeed, $V_1, \dots , V_m$ are $Beta(1,m-1)$-distributed and
$\Re(\lambda) >0$ guarantees that $\g E |V_1^{\lambda}| <\infty$.
Moreover, we are interested in solutions of (\ref{smoothing}) having a nonzero expectation and
the existence of such solutions implies that $\lambda$ is a root of (\ref{car}). Note that 
when $m\leq 13$, no $\lambda$ satisfies (\ref{hyp}).
\medskip

The following theorem implies Theorem~\ref{W}(i) because $W$ is integrable with
$\g EW=\frac{1}{\Gamma (1+\lambda_2 )}\neq 0$ (see~\cite{Pou05} for instance).

\begin{Th}
\label{support}
Let $\lambda$ be a non real complex number having a positive real part.
If $Z$ is a solution of~(\ref{smoothing}) having a nonzero expectation,
then the support of $Z$ is the whole complex plane.
\end{Th}

The proofs of Theorem~\ref{support} and Theorem \ref{density} make use of the complex-valued random variable
\begin{equation} \label{defA}
A:=\sum _{k=1}^mV_k^\lambda .
\end{equation}
Notice that 
the existence of an integrable solution $Z$ of~(\ref{smoothing}) such that $\g E (Z) \neq 0 $ implies 
that
\be
\label{conservation}
\g E(A)=1, 
\ee
which just means that $\lambda $ is a root of the polynomial (\ref{car}).
\medskip

\noindent {\sc Proof of Theorem \ref{support}.} For a complex valued random variable $X$, we denote its support by
$$  \Supp (X) = \{ x\in \C,  \; \forall \varepsilon >0, \ \g P ( |X-x| < \varepsilon ) >0  \}. $$
Let $Z$ be a solution of~(\ref{smoothing}) having a nonzero expectation. We first prove that $\forall a\in\g C,~\forall z\in\g C$,
\be
\label{implicationSupport}
\Big[ a\in\Supp (A) ~{\rm and}~z\in\Supp (Z)\Big] \Longrightarrow az\in\Supp (Z).
\ee
Indeed, let $\varepsilon >0$, $a\in\Supp(A)$ and $z\in\Supp (Z)$.
Let also $Z_1,\dots ,Z_m$ be i.i.d. copies of $Z$.
Then, with positive probability, $|A-a|\leq\varepsilon$ and $|Z_k-z|\leq\varepsilon$ for any $k$.
Therefore, with positive probability,
$$
\left|\sum _{k=1}^mV_k^{\lambda}Z_k-az\right|
=\left|\sum _{k=1}^mV_k^{\lambda}\left( Z_k-z\right) +z(A-a)\right|
\leq(m+\varepsilon)\varepsilon +|z|\varepsilon .
$$
The positive $\varepsilon$ being arbitrary, this shows that
$az\in\Supp\left(V_1^{\lambda}Z_1+\dots +V_m^{\lambda}Z_m \right)$ which implies that
$az\in\Supp (Z)$ because of~(\ref{smoothing}).

Let $z\in\Supp (Z)\setminus\{ 0\}$.
Such a $z$ exists because $\g E(Z)\neq 0$.
Iterating~(\ref{implicationSupport}), any complex number of the form $a_1\dots a_nz$
where $a_1,\dots ,a_n\in\Supp (A)$ belongs to $\Supp (Z)$.
Therefore, Lemmas~\ref{interieurSupport} and~\ref{monoid} below imply that $\Supp (Z)$ contains
$\g C\setminus\{ 0\}$ which suffices to conclude since the support of a probability measure is a closed
set.
\QED

\begin{Lem}
\label{interieurSupport} 

There exist $c,c'\in\g C\setminus \{0\}$ and respective open neighbourhoods $V$ and $V'$ of
$c$ and $c'$ such that $|c|>1$, $|c'|<1$ and $V\cup V'\subseteq\Supp (A)$.
\end{Lem}

\pff
Obviously,
$$
\Supp (A)=\left\{ \sum _{k=1}^mt_k^\lambda ,~0\leq t_k\leq 1,~\sum _{k=1}^mt_k=1\right\} .
$$
In particular, $\Supp (A)$ contains the set $f \left( [0,1]^2\right)$ (the image of $[0,1]^2$ by $f$),  where
$f$ is defined by
$$
\begin{array}{rccl}
f:&[0,1]^2&\to&\g C\\
&(s,t)&\mapsto&(st)^\lambda +\left( s(1-t)\right) ^\lambda +(1-s)^\lambda .
\end{array}
$$
We show that there exist $(s_c,t_c)$ and $(s_{c'},t_{c'})$ in $]0,1[^2$ such that $c:=f (s_c,t_c)$
and $c':=f (s_{c'},t_{c'})$ satisfy $|c|>1$, $0<|c'|<1$ and $f$ is a local
diffeomorphism in some respective neighbourhoods of $(s_c,t_c)$ and $(s_{c'},t_{c'})$, which implies
the result.

Let $\sigma$ and $\tau$ be respectively the real part and the imaginary part of $\lambda$. By assumption, $ 0< \sigma <1$. We assume that $\tau >0$; if not, replace $Z$ and $\lambda$ by their conjugates.
For any integer $k\geq 1$, denote
$$
u_k=\exp \left(-\frac{2k\pi}{\tau }\right)
{\rm ~and~~}
u'_k=\exp \left(\frac{\pi -2k\pi}{\tau }\right) .
$$
Then, $u_k$ and $u'_k$ are reals in $]0,1[$ that tend to $0$ as $k$ tends to infinity, and they satisfy
$$
u_k^\lambda =u_k^\sigma\in ]0,1[
{\rm ~and~~}
{u'_k}^\lambda =-{u'_k}^\sigma\in ]-1,0[ .
$$
Denote moreover
$$
\left\{
\begin{array}{l}
\displaystyle s_k:=u_k+u_k^2,\hskip 25pt t_k:=\frac 1{1+u_k}\\
\displaystyle s'_k:=u'_k+{u'_k}^2,\hskip 23pt t'_k:=\frac 1{1+u'_k}.
\end{array}
\right.
$$
As  $0<\sigma <1$, we have
$$
\left|f (s_k,t_k)\right|
=\left| u_k^\lambda +{u_k}^{2\lambda}+\left(1-u_k-{u_k}^2\right) ^\lambda\right|
=1+u_k^\sigma +O\left( u_k\right)
$$
and
$$
\left|f (s'_k,t'_k)\right|
=\left| {u'_k}^\lambda +{u'_k}^{2\lambda}+\left(1-u'_k-{u'_k}^2\right) ^\lambda\right|
=1-{u'_k}^\sigma +O\left( u'_k\right)
$$
when $k$ tends to infinity, so that $\left|f (s_k,t_k)\right| >1$ and $0<\left|f (s'_k,t'_k)\right| <1$
when $k$ is large enough.
It remains to show that $f$ is a local diffeomorphism in neighbourhoods of
$(s_k,t_k)$ and $(s'_k,t'_k)$.
Let $\Phi :[0,1]^2\to\g R^2$ be defined as
$$
\Phi (s,t)=\Big( \Re f (s,t),\Im f (s,t)\Big) .
$$
It suffices to show that the Jacobian of $\Phi$ at suitable $(s_c,t_c)$ and $(s_{c'},t_{c'})$ does not
vanish to show that
$f$ is a local diffeomorphism at these points.
This sufficient condition is equivalent to requiring that
$$
\frac{\partial f}{\partial s}\times\overline{\frac{\partial f}{\partial t}}
$$
is non real (the overline denotes the complex conjugacy).

For any $k\geq 1$, after computation, one gets
$$
\begin{array}{l}
\displaystyle
\frac{\partial f}{\partial s}(s_k,t_k)\times\overline{\frac{\partial f}{\partial t}}(s_k,t_k)
=|\lambda |^2\frac{1+u_k}{s_ku_k}
\Big[ u_k^\lambda\left( 1+u_k^\lambda\right)-s_k(1-s_k)^{\lambda -1}\Big]
\Big[ u_k^{\overline\lambda +1}-u_k^{2\overline\lambda}\Big]\\ \\
\displaystyle\hskip 120pt
=|\lambda |^2\frac{1+u_k}{s_ku_k}
\Big[u_k^{\sigma} + u_k^{2\sigma} - s_k(1-s_k)^{\lambda -1}\Big]
\Big[u_k^{\sigma + 1} - u_k^{2\sigma}\Big].
\end{array}
$$
The above number is real if and only if $(1-s_k)^\lambda\in\g R$, \emph{i.e.} if and only if
$\tau\log (1-s_k)\in\pi\g Z$.
Since $s_k\neq 0$ and $s_k$ tends to zero when $k$ tends to infinity,
$\tau\log (1-s_k)\notin\pi\g Z$ as soon as $k$ is large enough.
Therefore, taking $(s_c,t_c)=(s_k,t_k)$ for $k$ large enough suffices to get the result.
An argument of the same kind applied to the sequence $(s'_k,t'_k)_k$ leads to the result
on the existence of $(s_{c'},t_{c'})$.
\QED

\begin{Lem}
\label{monoid}
Let $V$ and $V'$ be respectively open neighbourhoods of $c \in \C$ and $c' \in \C$ with $|c|>1$ and $0<|c'|<1$, which do not contain $0$. 
Let
$$
M:=\left\{ v_1v_2\dots v_n,~n\geq 1,~v_1,v_2,\dots ,v_n\in V\cup V'\right\}.
$$
Then $M=\g C\setminus\{ 0\}$.
\end{Lem}

\pff
Let $\ell$ and $\ell '$ be complex numbers such that $\Re\ell >0$ and $\Re\ell '<0$.
Let $U$ and $U'$ be respectively open neighbourhoods of $\ell$ and $\ell '$.
Denote by $\rond M$ the additive submonoid of $\g C/2i\pi\g Z$ generated by $U\cup U'\mod 2i\pi$; it is the set of classes
$$
\rond M:=\{ u_1+u_2+\dots + u_n \mod 2i\pi ,~n\geq 1,~u_1,u_2,\dots ,u_n\in U\cup U'\}.
$$
We prove hereunder that $\rond M=\g C/2i\pi\g Z$.
Taking the exponential, this suffices to prove the lemma.  

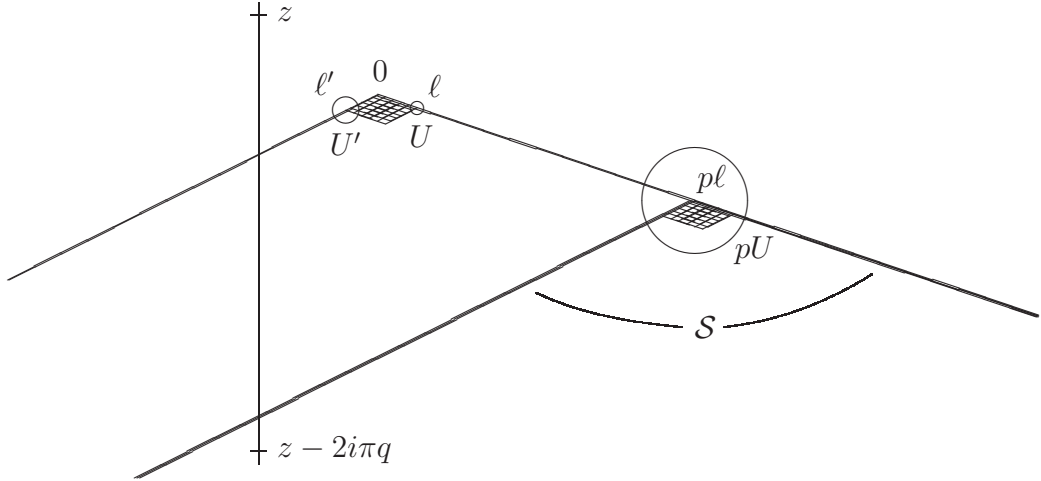
\begin{picture}(400,200)
\put(150,150){\line(-2,-1){140}}\put(150.5,149.75){\line(-2,-1){140}}
\put(150,150){\line(3,-1){250}}\put(150,149.5){\line(3,-1){250}}
\put(148,155){$0$}
\put(138,144){\line(3,-1){15}}\put(140,145){\line(3,-1){15}}\put(142,146){\line(3,-1){15}}\put(144,147){\line(3,-1){15}}\put(146,148){\line(3,-1){15}}\put(148,149){\line(3,-1){15}}
\put(165,145){\line(-2,-1){12}}\put(162,146){\line(-2,-1){12}}\put(159,147){\line(-2,-1){12}}\put(156,148){\line(-2,-1){12}}\put(153,149){\line(-2,-1){12}}
\put(165,145){\circle{5}}
\put(162,130){$U$}
\put(169,149){$\ell$}
\put(270,110){\circle{40}}
\put(271,115){$p\ell$}
\put(285,88){$pU$}
\put(258,104){\line(3,-1){15}}\put(260,105){\line(3,-1){15}}\put(262,106){\line(3,-1){15}}\put(264,107){\line(3,-1){15}}\put(266,108){\line(3,-1){15}}\put(268,109){\line(3,-1){15}}
\put(285,105){\line(-2,-1){12}}\put(282,106){\line(-2,-1){12}}\put(279,107){\line(-2,-1){12}}\put(276,108){\line(-2,-1){12}}\put(273,109){\line(-2,-1){12}}
\put(138,144){\circle{10}}
\put(132,127){$U'$}
\put(126,150){$\ell '$}
\put(270,110){\line(-2,-1){210}}\put(269,110){\line(-2,-1){210}}\put(268,110){\line(-2,-1){210}}
\put(270,109.6){\line(3,-1){130}}\put(269.2,109.3){\line(3,-1){130.6}}
\put(270,58){$\rond S$}
\qbezier(265,62)(230,65)(210,75)
\qbezier(281,62)(310,65)(337,82)
\put(112,178){$z$}
\put(105,185){\line(0,-1){175}}
\put(102,180){\line(1,0){6}}
\put(112,13){$z-2i\pi q$}
\put(102,15){\line(1,0){6}}
\end{picture}

Take an integer $p\geq 1$ large enough so that $pU$ contains a whole mesh of the lattice generated by
$\ell$ and $\ell'$, \emph{i.e.} such that $pU\supseteq p\ell+[0,1]\ell +[0,1]\ell '$.
Then, $\rond M$ contains the classes $\mod 2i\pi$ of the sector
$\rond S:=p\ell+\g R_{\geq 0}\ell+\g R_{\geq 0}\ell '$.
Since $\Re\ell\times\Re\ell '<0$, when $z$ is any complex number, there exists $q\in\g Z$
such that $z-2i\pi q\in\rond S$, which proves the result.
\QED

\section{Density and exponential moments}

As in the beginning of Section \ref{sec:support}, Theorem \ref{W}(ii) and (iii) are straightforward corollaries of the following theorem.
\begin{Th}
\label{density}
Let $\lambda$ be a non real complex number and $Z$ a solution of~(\ref{smoothing}) having a nonzero expectation.

\vskip 5pt
\hskip 10pt (i) If  $\Re(\lambda)>0$, then $Z$ admits a continuous square integrable density on $\g C$.

\vskip 5pt
\hskip 10pt (ii) If $\Re(\lambda)>\frac 12$, then
 $\g E e^{\delta  |Z|} < \infty$ for some $\delta >0$.
The exponential moment generating series of~$Z$ (thus) has a positive radius of convergence.
\end{Th}

\pff
It runs along the same lines as in \cite{ChaLiuPou} and uses the Fourier transform $\varphi$ of $Z$, namely
$$ \varphi(t) :=  \g E \exp \{ i \langle t,Z \rangle \}   = \g E \exp \{ i \Re (\overline t Z) \},    \quad t \in \C,    $$
where  $ \langle x,y\rangle   =   \Re (\overline x y ) = \Re (x) \Re (y)  + \Im (x) \Im (y) $.
In terms of Fourier transforms, Eq. (\ref{smoothing}) reads
\be
\label{equFourier}
 \varphi(t) =  \g E \left( \prod_{k=1}^m \varphi\left( tV_k^{\overline{\lambda}}\right) \right).
\ee
$\bullet$ To get (i), we prove that $\varphi$ is in $L^2(\g C)$  because it is dominated by $|t|^{-\delta}$ for some $\delta >1$ so that the inverse Fourier-Plancherel transform provides a square integrable density for $Z$. The guiding idea consists in adapting methods (developed in \cite{Liu99} and \cite{Liu01}) usually applied to positive real-valued random variables to the present complex-valued case.
For any $r\geq 0$, denote
\begin{equation*}
\psi(r) := \max_{|t| = r}|\varphi(t)|.
\end{equation*}
Using Theorem~\ref{support}, one can step by step mimick the proof of Theorem 7.17 in  \cite{ChaLiuPou} to get the result.
We just give hereunder an overview of this proof, written as successive hints.

Show first that Theorem~\ref{support} implies that $\psi (r)<1$ for any $r>0$.
Then, notice that
\be
\label{ineq}
\psi(r) \leq \g E \left( \prod_{k=1}^m \psi\left( r|V_k|^{\sigma}\right) \right).
\ee
By Fatou's lemma, (\ref{ineq})  implies that $\limsup_{+\infty}\psi (r)\in\{ 0,1\}$.
Iterating suitably inequality~(\ref{ineq}) leads to $\lim _{+\infty}\psi (r)=0$.
Finally, applying~(\ref{ineq}) again we can show  that $\psi(r) = O (r^{-\delta})$ for some $\delta >1$, so
  that $\varphi$ is square integrable on $\g C$, which
leads to the result.

\medskip

$\bullet$ To get (ii), like in \cite{ChaLiuPou}, we use Mandelbrot's cascades.
Denote $\rond V:= (V_{1}, V_{2}, \dots , V_{m})$.
Let $U$ be the set of finite sequences of positive integers between $1$ and $m$, namely
$$
U:= \bigcup_{n\geq 1} \{ 1, 2, \dots , m\}^n .
$$
Elements of $U$ are denoted by concatenation.
Let $\rond V_u:= (V_{u1}, V_{u2}, \dots , V_{um}), u\in U$ be independent copies of $\rond V$,
indexed by all finite sequences of integers $u= u_1\dots u_n \in U$.

Introduce the martingale $(Y_n)_{n\geq 1}$ defined by
$$
 Y_n :=  \sum_{u_1\dots u_{n} \in \{1,\dots , m\}^{n}}  V_{u_1}^{\lambda} V_{u_1u_2}^{\lambda} \dots V_{u_1\dots u_{n}}^{\lambda}.
$$

By (\ref{conservation}), $\g E(Y_n) = \g E(A) = 1$.
It can be easily seen that
\begin{equation}
\label{EqYn}
    Y_{n+1} = \sum_{k=1}^m V_k^{\lambda}Y_{n,k},
\end{equation}
where $ Y_{n,k} $ for $1\leq k\leq m $ are independent of each other and independent of the $V_k$ and  each has the same distribution as $Y_n$. Besides, since $\sigma > \frac 12$, $m\g EV_1^{2\sigma} < 1$ and, by Cauchy-Schwarz inequality,
$$
\g E |A|^2 \leq \g E\left( \sum_{k=1}^m |V_k^{\lambda}|\right)^2 = \g E\left( \sum_{k=1}^m V_k^{\sigma}\right)^2 \leq 2 \g E \sum_{k=1}^m V_k^{2\sigma} = 2m\g EV_1^{2\sigma} < 2.
$$
Therefore for $n\geq 1$, $Y_n$ is square integrable and
$$ \Var Y_{n+1} = (\g E |A|^2 - 1) + m\g EV_1^{2\sigma} \Var Y_n,  $$
where $\Var X = \g E\left( |X-\g EX|^2\right)$  denotes the variance of $X$.
Thus, the martingale  $(Y_n)_n$ is bounded  in $L^2$, so that  when $n\rightarrow + \infty$,
 \begin{equation*}
    Y_n \rightarrow  Y_\infty \mbox{ a.s. and in } L^2,
  \end{equation*}
where $Y_\infty$ is a (complex-valued) random variable with variance
 \begin{equation*}
  \Var (Y_\infty) = \frac {\g E |A|^2 - 1}{1  - m\g EV_1^{2\sigma}  }.
 \end{equation*}
Passing to the limit in Eq. (\ref{EqYn}) shows that $Y_\infty$ is a solution of Eq. (\ref{smoothing}) and by unicity, (ii) in Theorem \ref{density} holds as soon as it holds for  $Y_\infty$.

This last fact comes from an adaptation of Lemma 8.29 in \cite{ChaLiuPou}, giving some constants $C>0$ and $\varepsilon >0$ such that for
all $t \in \C$ with $|t| \leq \varepsilon $, we have
\begin{equation*}
  \g E  e^{\langle t, Y_\infty \rangle  } \leq e^{\Re (t)  + C |t|^2 }.
 \end{equation*}
The adaptation relies on $\sum_{k=1}^m V_k^{2\sigma} <1$ a.s. for $\sigma > \frac 12$. The last assertion implies 
that  $\g E  e^{t | Y_\infty | } < \infty$ for $t>0$ small enough,  so that the  exponential moment generating series of $Y_\infty $  has a positive radius of convergence.
\QED

\bibliographystyle{plain}
\bibliography{CLP}

\end{document}